\documentclass[a4paper,11pt]{amsart}
\usepackage{cases, cite}
\usepackage{amsfonts,latexsym,rawfonts,amsmath,amssymb,amsthm}
\usepackage[plainpages=false]{hyperref}
\usepackage{graphicx}
\usepackage{setspace}

\RequirePackage{color}

\theoremstyle{plain}
  \newtheorem{theorem}{Theorem}[section]
  \newtheorem{lemma}[theorem]{Lemma}
  \newtheorem{proposition}[theorem]{Proposition}
  
  \newtheorem{corollary}[theorem]{Corollary}

\theoremstyle{remark} 

\theoremstyle{definition}

\numberwithin{equation}{section}

\DeclareMathOperator*{\Ext}{Ext}
\DeclareMathOperator*{\supp}{supp}

\begin{document}

\title{Cubature formulas and Sobolev inequalities}

\author{Eli Putterman}
\address{School of Mathematical Sciences, Tel Aviv University, Tel Aviv, 66978, Israel.}
\email{putterman@mail.tau.ac.il}



\date{\today}



\begin{abstract}
We study a problem in the theory of cubature formulas on the sphere: given $\theta \in (0, 1)$, determine the infimum of $\|\nu\|_\theta = \sum_{i = 1}^n \nu_i^\theta$ over cubature formulas $\nu$ of strength $t$, where $\nu_i$ are the weights of the formula $\nu$. This problem, which generalizes the classical problem of bounding the minimal cardinality of a cubature formula -- the case $\theta = 0$ -- was introduced in recent work of Hang and Wang, who showed the problem to be related to optimal constants in Sobolev inequalities. Using the elementary theory of reproducing kernel Hilbert spaces on $S^{n - 1}$, we extend the best known upper and lower bounds for the minimal cardinality of strength-$t$ cubature formulas to bounds for the infimum of $\|\cdot\|_\theta$ for any $\theta \in (0, 1)$. In particular, we completely characterize the cubature measures of strength $3$ minimizing $\|\cdot\|_\theta$, showing that these are precisely the tight spherical $3$-designs.
\end{abstract}

\maketitle

\baselineskip=16.4pt
\parskip=3pt

\section{Introduction}
The approximation of integrals by finite sums, known as quadrature or cubature in higher dimensions, is a major branch of numerical analysis, and has connections to many areas in mathematics, such as representations of finite groups, orthogonal polynomials, lattices and sphere packings, and embeddings of Banach spaces; see \cite{HP} for a beautiful survey focusing on the case of Euclidean spheres. Recent work of Hang and Wang \cite{HW} has uncovered a connection between this subject and the theory of Sobolev inequalities, which raises a new set of questions in the theory of cubature formulas on the sphere.

Let $n \in \mathbb N$. For $p \in [1, n)$ we denote $p^* = \frac{np}{n - p}$. The classical Sobolev inequality in $\mathbb R^n$ states that there exists a constant $S_{n, p}$ such that
\begin{equation}
\|u\|_{p^*} \le S_{n, p} \|\nabla u\|_p
\end{equation}
for any $u \in W^{1, p}(\mathbb R^n)$. The optimal constant in this inequality is known to be
\begin{equation}\label{snp}
S_{n,p} = \frac{1}{n}\left( \frac{n\left( p-1\right) }{n-p}\right) ^{1-\frac{1}{p}}\left(\frac{n!}{\Gamma \left( \frac{n}{p}\right) \Gamma \left( n+1- \frac{n}{p}\right) \left\vert |S|^{n-1}\right\vert }\right) ^{\frac{1}{n}},
\end{equation}
and the functions saturating this inequality are explicitly known \cite{A1}.

Sobolev inequalities also hold for manifolds. On a compact Riemannian manifold $(M, g)$ of dimension $n$, the Sobolev inequality states that there exist constants $A_p, B_p$ (depending on $M$) such that for $u \in W^{1, p}(M)$, 
\begin{equation}
\|u\|_{p^*}^p \le A_p \|\nabla u\|_p^p + B_p \|u\|_p^p.
\end{equation}

The study of optimal constants in Sobolev and related inequalities has a long history; see \cite{He2} for an extensive treatment of Sobolev inequalities on manifolds, and see also \cite{BGL} for a heat-flow approach to many classical and recent functional inequalities. We mention only the results relevant to the current investigation. 

For compact manifolds, the optimal constants in Sobolev inequalities are basically known: Aubin \cite{A1} showed that for any $\varepsilon > 0$ there exists $c(\varepsilon)$ such that
\begin{equation}
\|u\|_{p^*}^p \le (S_{n, p}^p + \varepsilon) \|\nabla u\|_p^p + c(\varepsilon) \|u\|_p^p
\end{equation}
for all $u \in W^{1, p}(M)$. (In the case $p = 2$ one may take $\varepsilon = 0$, by a result of Hebey-Vaugon \cite{HV}.)

Under symmetry assumptions on $u$, the optimal constants in Sobolev inequalities can be improved. For instance, Hebey and Vaugon have shown \cite{He1} that if $M$ is acted upon by a compact group of isometries $G$ such that each point in $M$ has an orbit of cardinality at least $k$, then for any $\varepsilon$ there exists $c(\varepsilon)$ such that for any $G$-invariant function $u \in W^{1, p}(M)$,
\begin{equation}
\|u\|_{p^*}^p \le k^{-\frac{n}{p}} (S_{n, p}^p + \varepsilon) \|\nabla u\|_p^p + c(\varepsilon) \|u\|_p^p.
\end{equation}
For instance, even functions on $S^n$ satisfy such an improved Sobolev inequality with $k = 2$.

This can be generalized as follows: any even function $u \in W^{1, p}(S^n)$ has vanishing first moments, i.e., 
\begin{equation}
\int_{S^n} x_i u = 0
\end{equation}
for any $i = 1, \ldots, n + 1$. A theorem of Aubin \cite{A2} states that any function $u \in W^{1, p}(S^n)$ such that $|u|^{p^*}$ has vanishing first moments satisfies
\begin{equation}
\|u\|_{p^*}^p \le 2^{-\frac{n}{p}} (S_{n, p}^p + \varepsilon) \|\nabla u\|_p^p + c(\varepsilon) \|u\|_p^p.
\end{equation}
This result has applications to the Nirenberg problem of the existence of a metric on $S^n$ conformal to the usual one with prescribed scalar curvature.

Using the concentration-compactness principle of Lions \cite{Ln1,Ln2}, Hang and Wang \cite{HW} generalized the result of Aubin to higher moments. To state their result, we introduce some notation. 

For a nonnegative integer $k$, we denote%
\begin{eqnarray}
\mathcal{P}_{k} &=&\left\{ \text{all polynomials on }\mathbb{R}^{n+1}\text{
of degree at most }k\right\} ;  \label{eq1.10} \\
\overset{\circ }{\mathcal{P}}_{k} &=&\left\{ f\in \mathcal{P}_{k}:\int_{%
S^n}f \,d\sigma = 0\right\} .  \label{eq1.11}
\end{eqnarray}%
Here $\sigma$ denotes the standard (volume) measure on $S^n$. The restriction of $\mathcal P_k$ to $S^n$ consists of the linear span of spherical harmonics of order at most $k$, and $\overset{\circ}{\mathcal P}_k$ consists of the span of the spherical harmonics orthogonal to constant functions, or equivalently, the span of the spherical harmonics of order $1, \ldots, k$.

Let $\Delta(S^n)$ denote the space of discrete probability measures on $S^n$, i.e., the space of probability measures on $S^n$ supported on countably many points. For $t\in \mathbb N$, we denote%
\begin{eqnarray}
&&\mathcal{M}_{t}^{c}\left( S^n\right)  \label{eq1.12} \\
&=&\left\{ \nu :\nu \in \Delta(S^n) \,|\,\int_{S^n} f\,d\nu = 0 \text{ for all 
}f\in \overset{\circ }{\mathcal{P}}_t\right\}  \notag \\
&=&\left\{ \nu : \nu \in \Delta(S^n) \,|\,\frac{1}{|S^n|} \int_{S^n} f\,d\sigma = \int_{S^n} f\, d\nu 
\text{ for all }f\in \mathcal{P}_t\right\} .  \notag
\end{eqnarray}%
We call elements of $\mathcal M^c_t(S^n)$ spherical cubature formulas or cubature measures of strength $t$. (We use $t$ rather than $m$ as in \cite{HW} because this is the standard notation in the field of cubature formulas.) 

For $0 \le \theta < 1$, define a functional $\|\cdot\|_\theta: \Delta(S^n) \to [0, \infty]$ as follows:
\begin{equation}
\left\|\sum_i \nu_i \delta_{\xi_i} \right\|_\theta = \sum_i \nu_i^\theta.
\end{equation}
$\|\cdot\|_\theta$ is positive, nondegenerate, $\theta$-homogeneous, and strictly concave. In addition, $\|\nu\|_\theta$ is monotonically decreasing with increasing $\theta$, and $\|\nu\|_0$ is the cardinality of the support of $\nu$.

Hang and Wang define $\Theta(t, \theta, n)$ as the infimal value of $\|\nu\|_\theta$ for a spherical cubature formula $\nu$ of strength $t$ on $S^n$:
\begin{equation}\label{theta}
\Theta \left(t, \theta, n\right) = \inf \{\|\nu\|_\theta: \nu \in \mathcal M^c_t(S^n)\}
\end{equation}
As we shall see below (Corollary \ref{inf_att}), this infimum is in fact a minimum.

The main result of \cite{HW} states that the optimal constant in the Sobolev inequality on $S^n$ is improved by a factor of $\Theta(t, \frac{n - p}{n}, n)$ when restricting to $u$ such that $|u|^{p^*}$ is orthogonal to spherical harmonics of order up to $t$:

\begin{theorem}\label{hw_thm}
Assume $n\in \mathbb{N}$, $1 < p < n$ and $t \in \mathbb{N}$.
Denote $p^* = \frac{pn}{n-p}$. Then for any $\varepsilon > 0$, there exists $c(\varepsilon) > 0$ such that for any $u\in
W^{1,p}\left( S^n\right) $ with%
\begin{equation}\label{u_orth}
\int_{S^n} f |u|^{p^*}\, d\sigma = 0
\end{equation}%
for all $f\in \overset{\circ }{\mathcal{P}}_t$, we have%
\begin{equation}\label{impr_sob}
\|u\|_{p^*}^p \le \frac{S_{n, p}^p + \varepsilon}{\Theta\left(t, \frac{n-p}{n}, n\right)} \|\nabla u\|_p^p + c(\varepsilon) \|u\|_p^p, 
\end{equation}%
where $S_{n,p}$ is given by \eqref{snp} and $\Theta \left(t ,\frac{n-p}{n}, n\right) $ is given by \eqref{theta}.
\end{theorem}

Hang and Wang also proved, in the case $t = 2$, that $\Theta(t, \frac{n-p}{n}, n)$ is the best factor by which optimal constant in the Sobolev inequality can be improved for $u$ satisfying the conditions of the theorem, by constructing a family of such functions $u_\epsilon$ concentrated around the points of a cubature measure $\nu$ of strength $2$ attaining the bound $\|\nu\|_\theta = \Theta(2, \frac{n-p}{n}, n)$ which satisfy $\|u_\epsilon\|_{p^*} \ge \frac{S_{n, p}^p + \varepsilon}{\Theta(t, \frac{n-p}{n},n)} \|\nabla u_\epsilon\|_p + o(1)$ and $\|u_\epsilon\|_p = o(1)$. As their proof makes no use whatsoever of the structure of $\nu$ beyond the fact that it is a cubature measure, it goes through in general to show that $\Theta(t, \frac{n-p}{n}, n)$ is the best factor by which the constant in the Sobolev inequality can be improved for such $u$.

Hang and Wang showed that $\Theta(1, \theta, n) = 2^{1 - \theta}$, which recovers the result of Aubin, and proved that a cubature formula of strength $1$ minimizes $\|\cdot\|_\theta$ if and only if it is a uniform measure supported on a pair of antipodal points. They also proved that $\Theta(2, \theta, n) = (n + 2)^{1 - \theta}$, with minimizing measures given precisely by uniform measures on the vertices of a regular simplex inscribed in $S^n$, and that $\Theta(t, \theta, 1) = (t + 1)^{1 - \theta}$, with minimizing measures given by uniform measures on $t + 1$ evenly spaced points in $S^1$.


The present article continues the study of the constants $\Theta(t, \theta, n)$ initiated by \cite{HW}. We show that $\Theta(3, \theta, n) = (2n)^{1 - \theta}$ for any $n$ and explicitly describe the minimizers. We also give upper and lower bounds for $\Theta(t, \theta, n)$ in the regimes $n \ll t$, $t \ll n$, which match the upper and lower bounds known in the case $\theta = 0$. We summarize our results in the following theorem:

\begin{theorem}\label{main_results}\mbox{}
\begin{enumerate}
  \item For any $n$ and $\theta$, $\Theta(3, \theta, n) = (2n)^{1 - \theta}$, and $\nu \in \mathcal M^c_3(S^{n - 1})$ satisfies $\|\nu\|_\theta = (2n)^{1 - \theta}$ if and only if there exists an orthogonal basis $u_1, \ldots, u_n$ of $\mathbb R^n$ such that $\nu$ is the uniform measure on $\{\pm u_i\}_{i = 1}^n$.
  \item For any $t$, there exist constants $c, C$ such that 
  \begin{equation}
  (c n^{\lfloor \frac{t}{2} \rfloor})^{1 - \theta} \le \Theta(t, \theta, n) \le (C n^{2\lfloor\frac{t}{2} \rfloor - 1})^{1 - \theta}
  \end{equation}
   for all $n$ and $\theta$.
  \item For any $n$, there exist constants $c, C$ such that 
  \begin{equation}
  (c t^{n - 1})^{1 - \theta} \le \Theta(t, \theta, n) \le (C t^{n - 1})^{1 - \theta}
  \end{equation} for all $t$ and $\theta$.
\end{enumerate}
\end{theorem}

The main idea of the proof, in brief, is that the by now classical theory of Delsarte-Goethals-Seidel \cite{DGS}, which in its usual formulation gives bounds on the cardinality of cubature formulas -- a \textit{global} property of cubature formulas --  may be improved to a \textit{local} theory which bounds each individual weight of a cubature formula. In fact, this local formulation makes even the proofs of the classical results in the theory more transparent and intuitive, and also points the way to generalizations.

Following background material on spherical cubature in \S \ref{sec:prelim}, we prove the asymptotic upper bounds in Theorem \ref{main_results} in \S \ref{sec:first_prop}, give our version of the Delsarte-Goethals-Seidel theory and use it to prove the asymptotic lower bounds in \S \ref{sec:asymp}, and completely describe the minimizers in the case $t = 3$ in \S \ref{sec:3cub}. In \S \ref{sec:concl} we conclude by suggesting further directions for research.

\section{Preliminaries on spherical cubature}\label{sec:prelim}
We briefly survey the results about cubature formulas on the sphere that we shall need; all results not explicitly referenced may be found in \cite{HP} with full references and/or proofs. For compatibility with most of the literature on spherical cubature, we change notations from $S^n$ to $S^{n - 1}$, and use $d\sigma$ to denote the $O(n)$-invariant probability measure on $S^{n - 1}$ rather than the volume measure.

As above, we define a cubature measure on $S^{n - 1}$ of strength $t$ to be a discrete measure $\nu = \sum_{i = 0}^N \nu_i \delta_{\xi_i}$, where $N$ is finite or $\infty$, satisfying $\int_{S^{n - 1}} p\,d\nu = \int p\,d\sigma$ for all polynomials $p$ of degree at most $t$. Following \cite{HW} we write $\mathcal M^c_t(S^{n - 1})$ for the set of cubature measures of strength $t$. We use the terms ``cubature measure of strength $t$,'' ``$t$-cubature formula,'' and so forth interchangeably to refer to elements of $\mathcal M^c_t(S^{n - 1})$.

For example, a discrete measure $\nu$ is a cubature formula of strength $1$ if and only if its center of mass is the origin, and of strength $2$ if in addition its covariance matrix is proportional to the identity matrix (equivalently, its inertia ellipsoid is a ball); such a measure is also called an isotropic measure. The following characterization of cubature measures is fundamental:

\begin{theorem}\label{cub_equiv} Let $\nu = \sum_{i = 1}^N \nu_i \delta_{\xi_i}$ be a discrete probability measure on $S^{n - 1}$, and let $t \in \mathbb N$. The following are equivalent:
\begin{enumerate}
  \item $\nu$ is a cubature measure of strength $2t$.
  \item For any orthogonal transformation $T \in O(n)$ and any polynomial $p$ of degree at most $2t$, $\int_{S^{n - 1}} p\,d\nu = \int (p \circ T)\,d\nu$.
  \item $\int p_j\,d\nu = 0$ for any harmonic polynomial $p_j$ on $\mathbb R^n$ of order $j \in \{1, \ldots, 2t\}$.
  \item There exists a constant $c$ such that $\int_{S^{n - 1}} \langle x, u\rangle^{2t}\,d\nu(x) = c \langle u, u \rangle^t$ for any $u \in \mathbb R^n$, and $\int_{S^{n - 1}} \langle x, u\rangle^{2t - 1}\,d\nu(x) = 0$ for any $u \in \mathbb R^n$.
  \item $\int_{S^{n - 1}} \langle x, u\rangle^{2t}\,d\nu(x) = c_{2t, n}$ for any $u \in S^{n - 1}$, where 
  \begin{equation}\label{c_2t, n}
  c_{2t, n} = \frac{1\cdot 3 \cdot 5 \cdots (2t - 1)}{n(n + 2) \cdots (n + 2t - 2)}
  \end{equation}, and $\int_{S^{n - 1}} \langle x, u\rangle^{2t - 1}\,d\nu(x) = 0$ for any $u \in \mathbb R^n$.
\end{enumerate}

In addition, $\nu$ has strength $2t + 1$ if and only if it has strength $2t$ and satisfies $\int_{S^{n - 1}} p\,d\nu = 0$ for any homogeneous polynomial $p$ of degree $2t + 1$.
\end{theorem}

In particular, any even cubature measure $\nu$ of strength $2t$ is automatically of strength $2t + 1$; such $\nu$ are also called antipodal. Note that given any strength-$2t$ cubature formula $\nu$, we can construct an antipodal strength-$(2t + 1)$ cubature formula by taking the average of $\nu$ with its reflection through the origin.

Another way of putting the last two conditions in Theorem \ref{cub_num_bounds} is that if $\nu = \sum_{i = 1}^N \nu_i \delta_{\xi_i}$ is a spherical cubature formula of strength $2t$ then the map $u \mapsto (\left(\langle \frac{\nu_i}{c_{2t, n}}\right)^{\frac{1}{2t}}\xi_i, u\rangle)_{i = 1}^N$ is an isometric embedding of the Hilbert space $\ell^2_n$ into the Banach space $\ell^{2t}_N$. Conversely, given an isometric embedding $A: \ell^2_n \to \ell^{2t}_n$, the rows of $A$ yield a measure $\nu$ which satisfies $\int p\,d\nu = \int p\,d\sigma$ for all even polynomials of degree up to $2t$, which may be turned into a strength-$2t$ cubature formula by averaging it with its reflection through the origin. 

A central question in the theory of cubature formulas, with obvious applications to numerical integration, is the determination of cubature formulas of minimal cardinality. (We remark that the problem of computing $\Theta(t, \theta, n)$ can be thought of as a generalization of this question, since finding the minimum cardinality of a cubature formula of strength $t$ is the same as finding the minimum of $\|\cdot\|_0$ on $\mathcal M^c_t(S^n)$, i.e., computing $\Theta(t, 0, n)$.) The best fully general bounds known in the case of the sphere are given in the following theorem:

\begin{theorem}\label{cub_num_bounds} Let $\mathcal P_t^{(h)}(S^{n - 1})$ be the space of homogeneous polynomials of degree $t$ on $S^{n - 1}$, which has dimension $\binom{n + t - 1}{n - 1}$; we have $\mathcal P_t(S^{n - 1}) = \mathcal P_t^{(h)}(S^{n - 1}) \oplus \mathcal P_{t - 1}^{(h)}(S^{n - 1})$.

Let $M(t, n)$ denote the minimal number of points in a cubature formula of strength $t$. Then $\underline{N}(t, n) \le M(t, n) \le \overline{N}(t, n)$, where
\begin{align}
\underline{N}(2t, n) &= \dim P_t(S^{n - 1}) = \binom{n + t - 1}{n - 1} + \binom{n + t - 2}{n - 1}, \\
\overline{N}(2t, n) &= \dim \mathcal P_{2t}^{(h)}(S^{n - 1}) - 1 = \binom{n + 2t - 1}{n - 1} - 1,
\end{align}
and 
\begin{align}
\underline{N}(2t + 1, n) &= 2\dim \mathcal P_t^{(h)}(S^{n - 1}) = 2\binom{n + t - 1}{n - 1}, \\
\overline{N}(2t + 1, n) &= 2 \overline{N}(2t, n).
\end{align}
\end{theorem}

The lower bound for strength-$2t$ formulas, due to Stroud, comes from the fact that a cubature measure of order $2t$ defines a nondegenerate inner product $(p, q) \mapsto \int pq\,d\mu = \int pq\,d\sigma$ on $\mathcal P_t(S^n)$. The upper bound comes from a simple argument based on Fenchel's extension of Carathéodory's theorem stating that every point in the convex hull of a connected set $S \subset \mathbb R^N$ is a convex combination of $N$ points of $S$. (Similar ideas will be used in the proof of Proposition \ref{extr_pts} below.) The upper bound for the minimal cardinality of a strength $(2t + 1)$-formula follows immediately from the upper bound for strength $2t$-formulas, as any strength-$2t$ formula $\mu$ yields a strength-$(2t + 1)$ formula by averaging $\mu$ and its reflection through the origin; the lower bound, which is somewhat more subtle, is a result of M\"oller \cite{Mo}.

Note that for fixed $n$ and $t \to \infty$, $\underline{N}(t, n)$ and $\overline{N}(t, n)$ are both $O(t^{n - 1})$, where the implied constants depend on $n$, while for fixed $t$ and $n \to \infty$, $\underline{N}(t, n) = O(n^{\lfloor\frac{t}{2}\rfloor})$ and $\overline{N}(t, n) = O(n^{2\lfloor \frac{t}{2}\rfloor})$, where the constants depend on $t$.

An asymptotic improvement of the upper bound on $\overline{N}(t, n)$ for fixed $t$ and varying $n$ was obtained by Kuperberg \cite{Ku}:

\begin{theorem}\label{kup_bound} For any $t$ there exists a constant $C_t$ such that for any $n$ there exists a cubature formula with $C_t n^{2\lfloor \frac{t}{2} \rfloor - 1}$ points. 
\end{theorem}


\subsection{Spherical designs}
Spherical $t$-designs were introduced in \cite{DGS} in analogy to the notion of combinatorial designs over finite sets. A cubature measure $\nu$ of strength $t$ is called a spherical $t$-design if it is uniform on its support (in particular, the support of a $t$-design is finite). A spherical design of strength $t$ consisting of $\underline{N}(t, n)$ points is called tight. It is known that any strength-$2t$ cubature formula supported on $\underline{N}(2t, n)$ points is a spherical $2t$-design. We note also that the lower bounds for the minimal cardinality of cubature formulas in Theorem \ref{cub_num_bounds} were proven in the case of spherical $t$-designs in \cite{DGS}, apparently independently of work in the field of cubature formulas.

Only four infinite families of tight spherical designs are known. On $S^1$, $t + 1$ evenly spaced points form a tight spherical $t$-design. On $S^n$, a pair of antipodal points, the $n + 1$ vertices of a regular $n$-simplex, and the $2n$ vectors $\{\pm u_i\}_{i = 1}^n$ for an orthogonal basis $u_i$ are tight $1$-, $2$-, and $3$-spherical designs, respectively. Also, each of the families above are the only tight spherical designs in their respective classes of cubature formulas. Note that the first three families coincide with the measures proven by Hang and Wang to be the minimizers of $\Theta(t, \theta, 1)$, $\Theta(1, \theta, n)$, and $\Theta(2, \theta, n)$, respectively. 

Certain configurations of points with large symmetry groups, e.g., the vertices of regular polytopes and spherical shells of highly symmetric lattices such as the $E_8$ root lattice and the Leech lattice, are known to be tight $t$-designs on $S^{n - 1}$ for particular $t > 3$ and $n > 2$. Conversely, it it is known that tight spherical $t$-designs on $S^{n - 1}$, $n \ge 3$, exist only for $t \in \{1, 2, 3, 4, 5, 7, 11\}$; hence no $t$-cubature formula of order $\underline{N}(t, n)$, exists on $S^{n - 1}$ for $t$ not lying in this set.

The question of the existence of spherical $t$-designs of asymptotically optimal cardinality was the focus of much attention over several decades, culminating in a series of works of Bondarenko, Radchenko, and Viazovska \cite{BR,BRV} proving the following theorem:

\begin{theorem}For any $n \in \mathbb N$ there exists a constant $C_n$ such that for any $t$ there exists a spherical $t$-design with $C_n t^n$ points. 
\end{theorem}

See \cite{BB} for a survey on spherical designs and some generalizations from the perspective of algebraic combinatorics, and \cite{BG} for a more recent survey focusing on the relation of spherical designs to minimum-energy problems on spheres.

\section{First properties of \texorpdfstring{$\Theta(t, \theta, n)$}{Theta(t, theta, n)}}\label{sec:first_prop}
The space of cubature measures $\mathcal M^c_t(S^{n - 1})$ of strength $t$ on $S^{n - 1}$ is a convex subset of the space of Dirac probability measures on $S^{n - 1}$, and $\|\cdot\|_\theta$ is a strictly concave function on $\mathcal M^c_t(S^{n - 1})$. Hence, the infimum in the definition of $\Theta(t, \theta, n) = \inf_{\nu \in \mathcal M^c_t(S^{n - 1})} \|\nu\|_\theta$ (Equation \eqref{theta}) may be computed on the set of extreme points of $\mathcal M^c_t(S^{n - 1})$. Our first result concerns the structure of this set.

\begin{proposition}\label{extr_pts} Let $\nu = \sum_{i = 1}^N \nu_i \delta_{\xi_i} \in \mathcal M^{2t}_c(S^{n - 1})$, $\nu_i > 0$, and suppose $N > \overline{N}(t, n)$ (possibly infinite), where $\overline{N}(t, n)$ is defined in Theorem \ref{cub_num_bounds}. Then $\nu$ is not an extreme point of $\mathcal M^c_t(S^{n - 1})$.
\begin{proof}Let $p_1, \ldots, p_{\overline{N}(t, n)}$ be a basis for $\mathcal P_{2t}^{(h)}$. The linear map $A: \mathbb R^N \to \mathbb R^{\overline{N}(t, n) + 1}$ given by $(Av)_j = \sum_{i = 1}^N v_i p_j(\xi_i)$ has a nontrivial kernel; let $v$ be a nonzero vector in the kernel. Take $\epsilon > 0$ such that the measures $\nu_\pm = \sum_{i = 1}^N (\nu_i \pm \epsilon v_i) \delta_{\xi_i}$ have positive weights. Since $x \mapsto \langle x, u\rangle^{2t}$ is $2t$-homogeneous for every $u \in S^{n - 1}$, $\int \langle x, u\rangle^{2t}\,d\nu_{\pm}(x) = \int \langle x, u\rangle^{2t}\,d\mu$ for any $u$, so $\nu_{\pm}$ lie in $\mathcal M^c_t(S^n)$ by Theorem \ref{cub_num_bounds}. Hence $\nu = \frac{\nu_+ + \nu_-}{2}$ is not an extreme point of $\mathcal M^c_t(S^n)$.
\end{proof}
\end{proposition}

We can say a bit more: let $\mathcal M(S^{n - 1})$ be the space of Borel measures on $S^{n - 1}$ with the weak* topology. Then $\mathcal M^c_t(S^{n - 1})$ is contained in the compact subset $K$ of $\mathcal M(S^{n - 1})$ consisting of probability measures orthogonal to all polynomials $p \in \overset{\circ}{\mathcal P}_t$, and it is easily seen, by the same argument as in the proof of the proposition, that extreme points of $K$ must be discrete measures, hence $\Ext(K) = \Ext(\mathcal M(S^{n - 1}))$. The Krein-Milman theorem then yields that $K$ is the closure of the convex hull of $\Ext(\mathcal M(S^{n - 1}))$. 

\begin{corollary}\label{inf_att} The infimum in the definition of $\Theta(t, \theta, n)$ is achieved, i.e., there exists $\nu \in \mathcal M^c_t(S^{n - 1})$ with $\|\nu\|_\theta = \Theta(t, \theta, n)$.
\begin{proof}Let $N = \overline{N}(t, n)$. By the above, the infimum of $\Theta(t, \theta, n)$ on $\mathcal M^c_t(S^n)$ is the same as the infimum of  
$$\left\{\left.\sum_{i = 1}^N \nu_i^\theta: (\xi_i, \nu_i) \in (S^n \times [0, 1])^N\,\right|\,\sum_{i = 1}^N \nu_i p_j(\xi_i) = 0\,\,\forall p_j \in \overset{\circ}{\mathcal P}_m\right\},$$
which is attained by compactness.
\end{proof}
\end{corollary}

In particular, we can consider only finitely supported measures in the definition of $\Theta(t, \theta, n)$. This justifies our use of the term cubature, which properly refers only to approximation of integrals by finite sums.

We now proceed to our main topic, the determination of bounds for $\Theta(t, \theta, n)$. A first observation is that if $\nu$ is a spherical $t$-design of cardinality $N$, then $\|\nu\|_\theta = \sum_{i = 1}^N \frac{1}{N^\theta} = N^{1 - \theta}$; more generally, if $\nu = \sum_{i = 1}^N \nu_i \delta_{\xi_i}$ is any cubature formula of cardinality $N$, then 
\begin{equation}
\nu\|_\theta = \sum_{i = 1}^N \nu_i^\theta \ \le \left(\sum_{i = 1}^N (\nu_i^\theta)^{\frac{1}{\theta}}\right)^\theta \cdot \left(\sum_{i = 1}^N 1\right)^{1 - \theta} = N^{1 - \theta}
\end{equation} 
by H\"older's inequality. This means that existence results for cubature formulas of given cardinality (Theorem \ref{cub_num_bounds}, Theorem \ref{kup_bound}) translate immediately to upper bounds on $\Theta(t, \theta, n)$:

\begin{theorem}\label{theta_bds} \mbox{}
\begin{enumerate}
\item For any $n$, there exists a constant $C_n$ such that $\Theta(t, \theta, n) \le (C_n t^{n - 1})^{1 - \theta}$ for all $t$ and $\theta$. In fact, $\limsup_{t \to \infty} \frac{\Theta(t, \theta, n)}{t^{(n - 1)(1 - \theta)}} \le \left(\frac{2^{n - 1}}{(n - 1)!}\right)^{1 - \theta}$ for any $n$.
\item For any $t$, there exists a constant $C_t$ such that $\Theta(t, \theta, n) \le (C_t n^{2\lfloor \frac{t}{2}\rfloor})^{1 - \theta}$ for all $n$ and $\theta$. If we restrict to $n$ greater than $\frac{t}{c}$ for some constant $c$, we may take $C_t = \frac{(1 + 2 c)^t}{t!}$ for even $t$ and $C_n = 2 \frac{(1 + 2 c)^{t - 1}}{(t - 1)!}$ for odd $t$.  
\item For any $t$, there exists a constant $C_t'$ such that $\Theta(t, \theta, n) \le (C_t' n^{2\lfloor \frac{t}{2}\rfloor - 1})^{1 - \theta}$ for all $n$ and $\theta$.
\end{enumerate}
\end{theorem}

The quantitative estimates in parts (i) and (ii) of the theorem come from analyzing the asymptotics of $\binom{n + 2t - 1}{n - 1} = \binom{n + 2t - 1}{2t}$ for fixed $n$ and fixed $t$, respectively. Estimates for the constant $C_t'$ of part (iii) are less easy to obtain, since the method used in \cite{Ku} to construct the cubature formulas of order $O(n^{2\lfloor \frac{t}{2}\rfloor - 1})$ does not lend itself easily to making the constants explicit. The many explicit constructions of spherical designs and cubature formulas in the literature for specific $t$ and $n$ also provide upper bounds on $\Theta(t, \theta, n)$. 

However, lower bounds on the cardinality of spherical $t$-designs or more generally of $t$-cubature formulas, do not automatically give lower bounds on $\Theta(t, \theta, n)$, since the fact that a measure is supported on ``many'' points presents no obstacle to its having a small $\|\cdot\|_\theta$-norm. What is necessary to derive lower bounds on $\Theta(t, \theta, n)$ is to prove that any $t$-cubature formula is anti-concentrated on the sphere, i.e., a ``large'' fraction of the mass of such a formula is distributed among ``many'' points. Our method for doing so is to bound the maximum mass at any point in such a formula, which is the main goal of \S \ref{sec:asymp}. In the case $t = 3$, we can construct cubature formulas which attain the maximum possible weights and hence minimize $\|\cdot\|_\theta$ among $3$-cubature formulas, as we show in \S \ref{sec:3cub}.

\section{Reproducing kernels and the asymptotic lower bounds on \texorpdfstring{$\Theta(t, \theta, n)$}{Theta(t, theta, n)}}\label{sec:asymp}
The most natural way to prove our main result, Theorem \ref{weight_bounds}, on the weights of cubature formulas, is by using the formalism of reproducing kernel Hilbert spaces, which seems to have first been explicitly introduced into the subject by \cite{HP}. (We will remark at greater length on the history at the end of the section.) We give a brief and self-contained account of the elements of the theory of reproducing kernel Hilbert spaces that we shall need, which are easy to develop in the finite-dimensional case. For more complete accounts, see \cite{HP} and the literature cited there.

Let $\mathcal H \subset L^2(S^{n - 1})$ be some finite-dimensional real Hilbert space of continuous functions on the sphere equipped with the usual inner product $(f, g) \mapsto \int fg\,d\sigma$ (again, $\sigma$ is normalized so that $\sigma(S^{n - 1}) = 1$). We assume that $\mathcal H$ is preserved under orthogonal transformations, i.e., $f \circ U \in \mathcal H$ for every $f \in \mathcal H$, $U \in O(n)$. Let $f_1, \ldots, f_N$, $N = \dim V$, be an arbitrary orthonormal basis for $\mathcal H$. We define the reproducing kernel $K: S^{n - 1} \times S^{n - 1} \to \mathbb R$ of $\mathcal H$ as 
\begin{equation}\label{ker_def}
K(x, y) = \sum_{i = 1}^N f_i(x) f_i(y).
\end{equation}

A number of properties of $K(x, y)$ follow readily from this definition. First, for any $f \in \mathcal H$,
\begin{equation}\label{repr_ker}
\int K(x, y) f(y) \,d\sigma(y) =  \sum_{i = 1}^N \int f_i(y) f(y\,d\sigma(y) \cdot f_i(x) = \sum_{i = 1}^N \langle f, f_i\rangle f_i(x) = f(x),
\end{equation}
i.e., taking the inner product of a function $f$ with $K(\cdot, y)$ ``reproduces'' the value of $f$ at $y$, which is the reason for the terminology. By the Riesz representation theorem, this uniquely defines $K(\cdot, y)$ as an element of $\mathcal H$ for any $y \in S^{n - 1}$, which shows that $K(x, y)$ does not depend on the orthonormal basis used to define it.

In particular, $\int K(x, y)^2\,d\sigma(y) = K(x, x)$ for any $x \in S^{n - 1}$. Also, we have 
\begin{equation}\int K(x, x)\,d\sigma(x) = \sum_{i = 1}^N \int f_i(x)^2\,d\sigma(x) = N
\end{equation} 
because the $f_i$ are orthonormal.

Finally, we claim that $K(x, y)$ depends only on the distance between $x$ and $y$ on the sphere, or equivalently, $K(x, y) = K(Ux, Uy)$ for any $U \in O(n)$. Indeed, we have 
\begin{equation}K(Ux, Uy) = \sum_{i = 1}^N f_i(Ux) f_i(Uy) = \sum_{i = 1}^N (f_i \circ U)(x) (f_i \circ U)(y).
\end{equation} 
But $U$ is an isometry on $S^{n - 1}$ and hence precomposition with $U$ is an isometry on $L^2(S^{n - 1})$, so $(f_1 \circ U), \ldots, (f_n \circ U)$ is also an orthonormal basis of $\mathcal H$. Hence $(x, y) \mapsto K(Ux, Uy)$ satisfies the definition of a reproducing kernel on $\mathcal H$, so we must have $K(Ux, Uy) = K(x, y)$. In particular, $K(x, x)$ is independent of $x$, so $K(x, x) = \int K(x, x)\,d\sigma(x) = N$ for any $x \in S^{n - 1}$.

All the properties of $K$ but the last one are valid for functions on any space; only the last property of $K$ depends on the geometry of the sphere, but as we shall see, it is crucial for the subsequent developments.

A remark on notational conventions: since giving the distance between two points on the sphere is the same as giving their inner product, the last property of $K$ means that there is some continuous function $Q$ on $[-1, 1]$ such that $K(x, y) = Q(\langle x, y\rangle)$, with $Q(1) = N$. Usually the theory of reproducing kernels for spaces of spherical harmonics is formulated in terms of $Q$, which turns out to be a Gegenbauer polynomial.

The machinery of reproducing kernels allows us to easily prove the main theorem of the section:

\begin{theorem}\label{weight_bounds} Let $\nu$ be a cubature formula of strength $m$ on $S^{n - 1}$, and let $\underline{N}(m, n)$ be defined as in Theorem \ref{cub_num_bounds}. Then any weight $\nu_i$ of $\nu$ satisfies $\nu_i \le \underline{N}(m, n)^{-1}$. Furthermore, in the case $m = 2t + 1$ , if $\nu(u_i) = \underline{N}(m, n)^{-1}$ then also $\nu(-u_i) = \underline{N}(m, n)^{-1}$.
\begin{proof}The proof naturally splits into two cases, $m = 2t$ and $m = 2t + 1$. In the first case, let $\mathcal H = \mathcal P_t$, which has dimension $N = \underline{N}(2t, n)$, and let $K_y(x) = K(x, y)$ denote the reproducing kernel of $\mathcal H$, which is a polynomial of degree $t$ in the variable $x$. Given a point mass $\nu_i \delta_{u_i}$ contained in $\nu$, consider the integral $\int_{S^{n - 1}} K_{u_i}(x)^2\,d\nu$. Since $K_{u_i}(x)^2$ is a polynomial of degree $2t$, we have $\int_{S^{n - 1}} K_{u_i}(x)^2\,d\nu = \int_{S^{n - 1}} K_{u_i}(x)^2\,d\sigma$, which equals $K(u_i, u_i) = N$ by the reproducing kernel property. On the other hand, $\int_{S^{n - 1}} K_{u_i}(x)^2\,d\nu \ge \nu_i K_{u_i}(u_i)^2 = N^2\nu_i$ because $K(x, u_i)^2$ is nonnegative on $S^{n - 1}$, so we must have $\nu_i \le \frac{N}{N^2} = \frac{1}{N}$, as desired.

The case $m = 2t + 1$ is only slightly more involved. Let $\mathcal H = \mathcal P^{(h)}_t$ be the space of homogeneous polynomials of degree $t$, and let $K_y(x) = K(x, y)$ be its reproducing kernel as above. Let $N$ be the dimension of $\mathcal H$, which satisfies $N = \frac{1}{2} \underline{N}(2t + 1, n)$; our goal is to show that any point mass $\nu_i \delta_{u_i}$ contained in $\nu$ satisfies $\nu_i \le \frac{1}{2N}$. Let $H$ be the hemisphere $\{x \in S^{n - 1}: \langle x, u_i\rangle \le 0\}$, and consider the test functions $K_{u_i}(x)^2$ and $\langle x, u_i\rangle K_{u_i}(x)^2$ on $S^{n - 1}$. Since $K_{u_i}(x)$ is a homogeneous polynomial, $K_{u_i}(x)^2$ is necessarily even, implying that $\langle x, u_i\rangle K_{u_i}(x)^2$ is odd. Also, both $K_{u_i}(x)^2$ and $\langle x, u_i\rangle K_{u_i}(x)^2$ are of degree at most $2t + 1$, so we obtain the following two inequalities:
\begin{equation}\label{even_test}
N = \int_{S^{n - 1}} K_{u_i}(x)^2\,d\sigma = \int_{S^{n - 1}} K_{u_i}(x)^2\,d\mu \ge N^2 \nu_i + \int_H K_{u_i}(x)^2, 
\end{equation}

\begin{multline}\label{odd_test} 
0  = \int_{S^{n - 1}} \langle x, u_i\rangle K_{u_i}(x)^2\,d\sigma = \int_{S^{n - 1}} \langle x, u_i\rangle K_{u_i}(x)^2\,d\mu \\ 
\ge N^2 \nu_i + \int_H \langle x, u_i\rangle K_{u_i}(x)^2\,d\nu \ge N^2 \nu_i - \int_H K_{u_i}(x)^2\,d\nu 
\end{multline}
where in the last inequality we have used the fact that $|\langle x, u_i\rangle| K(x, u_i)^2 \le K(x, u_i)^2$ on $S^{n - 1}$. Summing these two inequalities we obtain $2 N^2 \nu_i \le N$, or $\nu_i \le \frac{1}{2N}$, as desired. 

Finally, if $\nu_i = \frac{1}{2N}$ then equality must hold throughout \eqref{odd_test}, which, by examining the last estimate, implies that $\int_H \langle x, u_i\rangle K_{u_i}(x)^2\,d\nu = -\int_H K_{u_i}(x)^2\,d\nu$. This means that for every $u \in H$ such that $\nu(u) > 0$ we must have either $K_{u_i}(u)^2 = 0$ or $\langle u, u_i\rangle = -1$. We must also have equality between the extreme left- and right-hand sides of \eqref{odd_test}, i.e., $\int_H K_{u_i}(x)^2\,d\nu = N^2 \nu_i = \frac{N}{2}$. This implies that the weight of $\nu$ at the unique point $u = -u_i$ for which $\langle u, u_i\rangle = -1$ must be $\frac{1}{K_{u_i}(-u_i)^2} \cdot \frac{N}{2} = \frac{1}{2N}$ (again, using the fact that $K_{u_i}$ is even).
\end{proof}
\end{theorem}

The theorem implies the bounds of Theorem \ref{cub_num_bounds} on the minimal cardinality of cubature formulas as an immediate corollary; moreover, one obtains that whether $t$ is even or odd, the only strength-$t$ cubature formulas of minimal cardinality are tight spherical designs, which seems to have first been proven in \cite{LV}; moreover, for $t$ odd, one also obtains that such designs must be antipodal, which was proven in \cite{DGS}.

The asymptotic lower bounds on $\Theta(t, \theta, n)$ in Theorem \ref{main_results} follow immediately from the asymptotics of $\underline{N}(t, n)$ stated below Theorem \ref{cub_num_bounds}, with the help of a simple lemma:

\begin{lemma}\label{weight_bd}
If $\nu = \sum_{i = 1}^N \nu_i \delta_{\xi_i}$ is a probability measure with $\nu_i \le \alpha^{-1}$ for all $i$, then $\|\nu\|_\theta \le \alpha^{1 - \theta}$, with equality if and only if $\nu_i = \alpha^{-1}$ for all $i$.
\begin{proof}$\|\nu\|_\theta = \sum_{i = 1}^N \nu_i^\theta \ge \sum_{i = 1}^N \nu_i \cdot \alpha^{1 - \theta} = \alpha^{1 - \theta}$.
\end{proof}
\end{lemma}

Before moving on, we briefly comment on the relationship of the proof given here to preceding treatments. The use of Gegenbauer polynomials, which are reproducing kernels on spaces of spherical harmonics, goes back to Delsarte, Goethals, and Seidel \cite{DGS}, who derived the lower bounds on the cardinality of spherical $t$-designs for both odd and even $t$ by considering matrices obtained by evaluating a basis of spherical harmonics at points of the design; Gegenbauer polynomials appear when multiplying these matrices together. Further developments along these lines can be found in \cite{KL}, who discussed kernels in greater generality and on spaces other than $S^{n - 1}$, and in \cite{LV}, which introduced arguments based on kernels to the more general context of cubature formulas (cf. especially \cite[Theorem 4.5]{LV}, where a kernel makes an appearance in the argument that any tight $2t$-cubature formula is a $2t$-design), and in \cite{HP}, where the terminology of the theory of reproducing kernels is used for the first time. However, even in \cite{HP} the reproducing kernels are not used to bound individual weights of cubature formulas, but only the minimal cardinality, and furthermore, the argument in \cite{HP} is only for even strengths.

It seems that since the goal in these works was to bound the cardinality of strength-$t$ formulas or designs, which is a ``global'' property, it was not noticed that kernels provide ``local'' bounds on individual weights. The main contribution of this work, which is motivated by the fact that to obtain lower bounds on $\Theta(t, \theta, n)$ for $\theta > 0$ one needs bounds on the weights, is to emphasize that the control that theory of reproducing kernels gives over cubature formulas is actually more precise than hitherto observed; even the result that tight spherical $(2t + 1)$-designs must be antipodal turns out to be a consequence of a ``local'' result. 

We give more details on the relationship of our arguments to those of \cite{DGS} at the very end of \S \ref{sec:concl} 

\section{Minimizing cubature formulas of strength \texorpdfstring{$3$}{3}}\label{sec:3cub}
The main result of this brief section identifies the minimizers of $\Theta(t, \theta, n)$ in the case $t = 3$ as the tight spherical $3$-designs, complementing the results of Hang and Wang for $t = 1, 2$ and for $n = 1$. Thus, all the infinite families of tight spherical designs known are also the minimizers of the $\|\cdot\|_\theta$-norms in their respective $\mathcal M^c_t(S^n)$.

Before treating the case $t = 3$, we remark that Lemma \ref{weight_bd} gives a substantially simpler proof that $\Theta(1, \theta, n) = 2^{1 - \theta}$ with minimizers given only by measures of the form $\nu_\xi = \frac{1}{2} (\delta_\xi + \delta_{-\xi})$. Indeed, if $\|\nu\|_\theta \le 2^{1 - \theta}$ then we must have $\nu(\xi) \ge \frac{1}{2}$ for some $\xi$ by Lemma \ref{weight_bd}, and by the last statement of Theorem \ref{weight_bounds}, this implies that $\nu(-\xi) = \frac{1}{2}$, i.e., $\nu = \nu_\xi$. A similar argument will be used to treat the case of $3$-cubature formulas.

\begin{theorem}\label{cub3_opt}
For any $n \in \mathbb N$ and $\theta \in (0, 1)$, $\Theta(3, \theta, n) = (2n)^{1 - \theta}$, and $\nu \in \mathcal M^c_3(S^{n - 1})$ is a minimizing measure if and only if it is a tight spherical $3$-design, i.e., $\nu = \nu_{\{u_i\}}$ is the uniform measure supported on the union of an orthogonal basis $\{u_i\}_{i = 1}^n$ and its reflection through the origin.
\begin{proof}Taking $m = 3$ in Theorem \ref{weight_bounds} shows that any weight of any $\nu \in \mathcal M^c_3(S^{n - 1})$ is at most $\frac{1}{2n}$, so Lemma \ref{weight_bd} gives $\Theta(3, \theta, n) \ge (2n)^{1 - \theta}$, with equality if and only if all the weights of $\nu$ are $\frac{1}{2n}$. For any orthogonal basis $\{u_i\}$, $\nu_{\{u_i\}}$ is a $3$-cubature formula with $\|\nu_{\{u_i\}}\| = (2n)^{1 - \theta}$, so all that remains is to show that these are the only minimizers.

Let $\nu$ be a $3$-cubature formula minimizing $\|\nu\|_\theta$; by the above, $\nu$ must be the uniform measure on $2n$ points. Let $u$ such that $\nu(u) = \frac{1}{2n}$; without loss of generality we take $u = e_n$. By the final statement in Theorem \ref{weight_bounds} we have that $\nu(-e_n) = \frac{1}{2n}$ as well. Also, by Theorem \ref{cub_equiv}, we have $\frac{1}{n} = \int_{S^{n - 1}} x_n^2\,d\nu \ge \nu(e_n) + \nu(-e_n) = \frac{1}{n}$, which means that we must have 
\begin{equation}
\supp \mu \backslash \{-e_n, e_n\} \subset S^{n - 2} = \{x \in S^{n - 1}: x_n = 0\}.
\end{equation}

Consider the probability measure $\nu'$ obtained from $\nu$ by removing the mass at $\pm e_n$ and rescaling, i.e., $\nu' = \frac{n}{n - 1} \left(\nu - \frac{1}{2n} (\delta_{e_n} + \delta_{-e_n})\right)$. It's easy to see that $\nu'$ is a cubature formula of strength $3$ on $S^{n - 2}$ satisfying 
\begin{align}\|\nu'\|_\theta &= \left(\frac{n}{n - 1}\right)^\theta \left(\|\nu\|_\theta - 2\cdot \frac{1}{(2n)^\theta}\right) \nonumber \\
&\le \left(\frac{n}{n - 1}\right)^\theta ((2n)^{1 - \theta} - 2 (2n)^{-\theta}) = (2(n - 1))^{1 - \theta}.
\end{align}
Hence, by induction, $\nu'$ must be the uniform probability measure supported on $\{\pm v_i\}_{i = 1}^{n - 1}$ for $\{v_1, \ldots, v_{n - 1}\}$ an orthogonal basis of $\mathbb R^{n - 1} \subset \mathbb R^n$, so $\nu = \frac{1}{2n} (\delta_{e_n} + \delta_{-e_n}) + \frac{n - 1}{n} \nu'$ is a uniform probability measure supported on $\{\pm u_i\}_{i = 1}^n$ where $\{u_1, \ldots, u_n\} = \{v_1, \ldots, v_{n - 1}\} \cup \{e_n\}$ is an orthogonal basis of $\mathbb R^n$, as desired.
\end{proof}
\end{theorem}

The case $\theta = 0$ of the theorem is that any cubature formula of strength $3$ having cardinality $2n$ is the uniform measure supported on an orthogonal basis and its reflection through the origin; this was first proven by Mysovskikh \cite{My}.

\section{Conclusion and further directions}\label{sec:concl}
Following the pioneering work of Hang and Wang \cite{HW}, in this work we have obtained an essentially complete description, given the current state of knowledge regarding cubature formulas on the sphere, of the quantities $\Theta(t, \theta, n)$, which control the improvement in Sobolev inequalities on the sphere for functions orthogonal to spherical harmonics up to a given order. To conclude, we will mention some possibilities for generalizing these results.

First, it is entirely straightforward to generalize our arguments to the setting of compact two-point homogeneous manifolds (see \cite{W}), which are precisely those manifolds for which the reproducing kernel $K(x, y)$ can be written as $K(d(x, y))$, where $d$ is the metric. This setting was already introduced by \cite{KL} in order to generalize the theory in \cite{DGS}, and cubature formulas on compact two-point homogeneous manifolds were studied by \cite{BD} in order to obtain bounds on Kolmogorov $n$-widths of Sobolev spaces on these manifolds. 

For instance, using the fact that spaces of polynomials on the real projective spaces $\mathbb P^{n - 1}(\mathbb R)$ are precisely the spaces of even-degree polynomials on $S^{n - 1}$, the argument of Theorem \ref{weight_bounds} shows that the weight of any point in a cubature formula of strength $4t$ on $\mathbb P^{n - 1}(\mathbb R)$ is at most $(\dim \mathcal P_{2t}^{(h)})^{-1}$, where $\mathcal P_{2t}^{(h)}$ is the space of even-degree polynomials on $S^{n - 1}$ of degree at most $2t$, with a corresponding improvement in the optimal constants in the Sobolev inequality on $\mathbb P^{n - 1}(\mathbb R)$ for functions orthogonal to polynomials of degree $4t$.

In fact, our arguments do not rely on the two-point property, so they generalize to any compact space $X$ with a transitive group of isometries, given a sequence finite-dimensional subspaces of functions $\mathcal F_t \subset L^2(X)$ satisfying 
\begin{equation}\label{poly_filt}
\mathbb R = \mathcal F_0 \subset \mathcal F_1 \subset \mathcal F_2 \cdots
\end{equation}
and 
\begin{equation}\label{poly_mult}
\mathcal F_t \cdot \mathcal F_t \subset \mathcal F_{2t}.
\end{equation}
If $X$ is a homogeneous space of a compact Lie group $G$, the $\mathcal F_t$ may be constructed as direct sums of certain irreducible representations of $G$ in the Peter-Weyl decomposition of $L^2(X)$, using the theory of tensor products of representations to show that the condition \eqref{poly_mult} holds. It may be more natural, however, to abandon the ``one-dimensional'' filtration \eqref{poly_filt} entirely and develop the theory of cubature formulas with respect to the decomposition of $L^2(X)$ into unitary representations of $G$ directly.

A further generalization is to the case when the measure on $S^{n - 1}$ is not the homogeneous measure $d\sigma$, but rather $w\,d\sigma$ for some sufficiently nice weight function $w$. Cubature formulas $(S^{n - 1}, w\,d\sigma)$ were studied by Dai \cite{D} (see also the book of Dai and Xu \cite{DX}), who obtained bounds on the mass of cubature formulas in small balls which are asymptotically of the same order as the bounds on weights proved in this note, under the condition that $w$ is a doubling weight. The method used in \cite{D} to obtain such bounds is actually quite similar to the method employed here: one integrates a $2t$-cubature formula for $w\,d\sigma$ against the square of the usual kernel of $\mathcal P_t$, using the fact that $w$ is a doubling weight to control the integrals in terms of integrals with respect to the usual measure. We can thus improve Sobolev inequalities for functions on the sphere whose $p^*$-th powers are orthogonal to low-degree polynomials with respect to well-behaved measures other than area measure, and the factor by which the optimal constant is improved will be of the same asymptotic order as $\Theta(t, \frac{n - p}{n}, n)$.

Several more challenging directions for further research also suggest themselves. First, is it possible to obtain interesting bounds on cubature formulas for general Riemannian manifolds with respect to suitable finite-dimensional spaces of functions? For the method used in this work to apply, one would have to be able to bound the diagonal elements of reproducing kernels for such spaces, and it is not immediately clear how to do so in the non-homogeneous setting where $K(x, x)$ may depend on $x$.

Secondly, it is possible to generalize the main result of \cite{HW} to situations in which the function $u$ is not known to satisfy $\int f\,|u|^{p^*} = 0$ for all $f \in \mathcal P_t$ orthogonal to constants, but merely 
\begin{equation}\label{ineq_constr}
\int f\,|u|^{p^*} \le b(f) \|u\|_{p^*}^{p^*}
\end{equation} 
for some $b(f)$. (Of course, this is always satisfied with $b(f) = \|f\|_\infty$, but we are interested in situations where a better bound is known.) In this case, the method of \cite{HW} yields that the optimal constant in the Sobolev inequality can be improved for such $u$ by the infimum of $\|\nu\|_{\frac{n - p}{n}}$ over all discrete probability measures $\nu$ such that $\int f\,d\nu\le b(f)$ for all $f \in \overset{\circ}{\mathcal P_t}$. 

We give an example of a case in which spaces of ``almost-cubature measures'' arise naturally. Dvoretsky's theorem in asymptotic geometric analysis implies that for any $\epsilon$ there exists an embedding of $\ell^2_n$ into $\ell^{2t}_{c(\epsilon) n^t}$ with distortion at most $\epsilon$ for some function $c(\epsilon)$ independent of $n$ (see, e.g., \cite[Chapter 5]{AGM}). (Note that the asymptotics of $c(\epsilon) n^t$ precisely match the asymptotics of the lower bound for the minimal cardinality of a cubature formula of strength $t$, which corresponds to an exact embedding.) The condition that the embedding $\ell^2_n \hookrightarrow \ell^{2t}_{c(\epsilon) n^t}$ have distortion at most $\epsilon$ translates to a certain property of the measure $\nu$ constructed from the rows of the matrix of the embedding, as in the case where $\nu$ is a bona fide cubature formula: for polynomials of the form $p_u(x) = \langle x, u\rangle^{2t}$, $\left|\int p_u\,d\nu - \int p_u\,d\sigma\right| \le \epsilon \int p_u\,d\sigma$. Since these polynomials span $\mathcal P_{2t}^{(h)}$, one sees that the space of Dvoretsky-type embeddings $\ell^2_n \hookrightarrow \ell^{2t}_{c(\epsilon) n^t}$ corresponds to a space of discrete measures satisfying
\begin{equation}
\left|\int f\,d\nu - \int f\,d\sigma\right| \le \epsilon \cdot c(f)
\end{equation} 
for all $f \in \mathcal P_t^{(h)}$, where $c(f)$ can be estimated by the sum of the coefficients in an expansion $f = \sum_i c_i p_{u_i}$ in elements of the spanning set $\{p_u: u \in S^{n - 1}\}$. Unfortunately, we do not have good control on the coefficients of such an expansion, so $c(f)$ may be much larger than the norm of $f$, unlike $c(p_u) = \int p_u\,d\sigma$. Using the reproducing kernel $K_x(y) = K(x, y)$ of $\mathcal P_t^{(h)}$ and its defining property $f(x) = \int K_x(y) f(y)\,dy$, we can at least reduce the problem of bounding $c(f)$ to the problem of bounding $c(K_x)$. However, dimensional considerations suggest that $c(K_x)$ grows rapidly with $n$ and $t$, so for fixed $\epsilon$ the bound $\left|\int K_x\,d\nu - \int K_x\,d\sigma\right| \le \epsilon \cdot c(K_x)$ is no better than the trivial bound, which has $\|K_x\|_\infty$ on the right-hand side, except possibly for small $n$ and $t$. In short, it does not seem that Dvoretsky-type embeddings correspond to a natural set of constraints on functions in $W^{1, p}(S^{n - 1})$, for which it would be interesting to find the improvements in optimal constants in Sobolev inequalities for functions satisfying these constraints. But there may be other spaces of functions defined by inequality constraints of the form \ref{ineq_constr} for which the method proposed in \cite{HW} and developed here becomes relevant.

Finally, we may ask whether the bounds of Theorem \ref{weight_bounds} on the weights of cubature formulas can be improved. Essentially, the reason that integrating a reproducing kernel against a cubature formula $\nu$ gives good bounds on the weights of $\nu$ is that the kernel $K_u(x)$ is highly localized around $u$; indeed, $\frac{K_u(u)^2}{\|K_u\|_2^2} = N$, where $N$ is the dimension of the Hilbert space. But it is not necessarily the case that $K_u^2$ is the ``best'' test function, which would be the polynomial in $\mathcal P_{2t}(S^{n - 1})$ satisfying $p \ge 0$ and having $\frac{p(u)}{\int p}$ as large as possible. 

This question is related to the ``linear programming method'' introduced, in the context of spherical codes and designs, by \cite{DGS} and exploited spectacularly by \cite{KL} to provide upper bounds for sphere packings. In our context, it has been known since \cite{DGS} that any nonnegative function $F: [-1, 1] \to \mathbb R^+$ satisfying
\begin{equation}\label{test_fn_cond}
F = \sum_{k = 0}^\infty a_k P_k^{(\frac{n - 2}{2})}, \qquad \text{$a_k < 0$ for $k > t$,}
\end{equation}
where $P_k^{(\frac{n - 2}{2})}$ are Gegenbauer polynomials (namely, reproducing kernels of spherical harmonics of degree $k$), gives a bound of $\frac{F(1)}{a_0}$ for the minimum cardinality of a spherical $t$-design. ($a_0$ can be computed as $\int F\,d\mu_n$, where $d\mu_n = (1 - x^2)^{\frac{n - 3}{2}}\,dx$, the Gegenbauer weight function, is the one-dimensional projection of the normalized area measure on $S^{n - 1}$.) The reason for this is that any reproducing kernel $K$ on a finite-dimensional Hilbert space of functions $\mathcal H \subset S^{n - 1}$ satisfies
\begin{equation}\label{pos_def}
\int K(x, y) \,d\nu(x) \,d\nu(y) = \left(\sum_{i = 1}^N \int f_i(x)\,d\nu(x)\right) \left(\sum_{i = 1}^N \int f_i(y)\,d\nu(y)\right) \ge 0
\end{equation}
for $f_i$ an orthonormal basis of $\mathcal H$, since $K(x, y) = \sum f_i(x) f_i(y)$ (this property of reproducing kernels is known as positive-definiteness). Decompose $F(\langle x, y\rangle)$ as $\sum a_k Q_k(x, y)$ where the $Q_k$ are reproducing kernels of the spaces of spherical harmonics of degree $k$; in particular, $Q_0 = 1$ and so $a_0 = \int F\,d\mu_d$. Now let $\nu = \sum_{i = 1}^N \nu_i \delta_{\xi_i}$ be a cubature formula of strength $t$, and integrate the two-variable test function $F(\langle x, y\rangle)$ against $\nu \otimes \nu$:
\begin{align}
F(1) &= F(1) \sum_{i = 1}^N \nu_i^2 \le \int F(\langle x, y\rangle)\,d\nu(x)\,d\nu(y) = \sum_{k = 0}^\infty a_k \int Q_k(x, y)\,d\nu(x)\,d\nu(y) \nonumber \\
&= a_0 + \sum_{k > t} a_k \int Q_k(x, y) \,d\nu(x) \,d\nu(y) \le a_0,
\end{align}
where we have used the fact that $\int Q_k(x, y)\,d\nu(x)\,d\nu(y) = 0$ for $0 < k \le t$ because for any $y$, $Q_k(\cdot, y)$ is a polynomial of degree $k \le t$ and vanishing integral over the sphere, along with the positive-definiteness of each $Q_k$. This shows that $\sum_{i = 1}^N \nu_i^2 \le \frac{a_0}{F(1)}$; since $\sum_i \nu_i = 1$, Cauchy-Schwarz yields $\frac{1}{N} \le \sum_{i = 1}^N \nu_i^2 \le \frac{a_0}{F(1)}$, which gives the desired lower bound on the cardinality of strength-$t$ cubature formulas (and $t$-designs in particular).

The point here is that this argument does \textit{not} ``localize'' in the spirit of Theorem \ref{weight_bounds} to prove that each individual weight of the cubature formula is bounded above by $\frac{1}{F(1)}\int F\,d\mu_n$, since the positive-definiteness property of reproducing kernels has no analogue for single integrals of the form $\int_{S^{n - 1}} Q_k(x, y)\,d\nu(y)$. The method of proof of Theorem \ref{weight_bounds} thus requires a test function whose Gegenbauer coefficients actually vanish for $k > t$, which dictates the choice of the test functions used in the proof. For fixed $t$ and increasing $n$, Yudin \cite{Y} constructed test functions $F$ satisfying the conditions of \eqref{test_fn_cond} which are asymptotically more localized then the functions $K_u^2$ considered in the proof of Theorem \ref{weight_bounds}, and hence give better bounds on the cardinality of strength-$t$ cubature formulas. However, we do not know if $K_u^2$ is the best choice of test function among polynomials of degree at most $2t$; if not, then we could improve Theorem \ref{weight_bounds}. In addition, perhaps a modification of the proof of Theorem \ref{weight_bounds} would enable it to use test functions $F$ whose Gegenbauer components of order $k > t$ do not necessarily vanish; this might involve a detailed analysis of the behavior of the higher-order components of $F$ at the points of the formula. We leave the exploration of such possibilities to future work.

\end{document}